\documentclass[twoside,11pt]{article}
\usepackage{geometry}
\usepackage{graphicx}
\usepackage{indentfirst}                                
\usepackage[colorlinks]{hyperref}
\hypersetup{linkcolor=blue,filecolor=black,urlcolor=blue, citecolor=black}   
\usepackage[numbers,sort&compress]{natbib}         

\ifxetex
\usepackage{letltxmacro}
\LetLtxMacro\SavedIncludeGraphics\includegraphics
\def\includegraphics#1#{
 \IncludeGraphicsAux{#1}%
}%
\newcommand*{\IncludeGraphicsAux}[2]{%
 \XeTeXLinkBox{%
  \SavedIncludeGraphics#1{#2}%
}}
\fi
\newcommand\orcidicon[1]{\href{https://orcid.org/#1}{\includegraphics[scale=0.02]{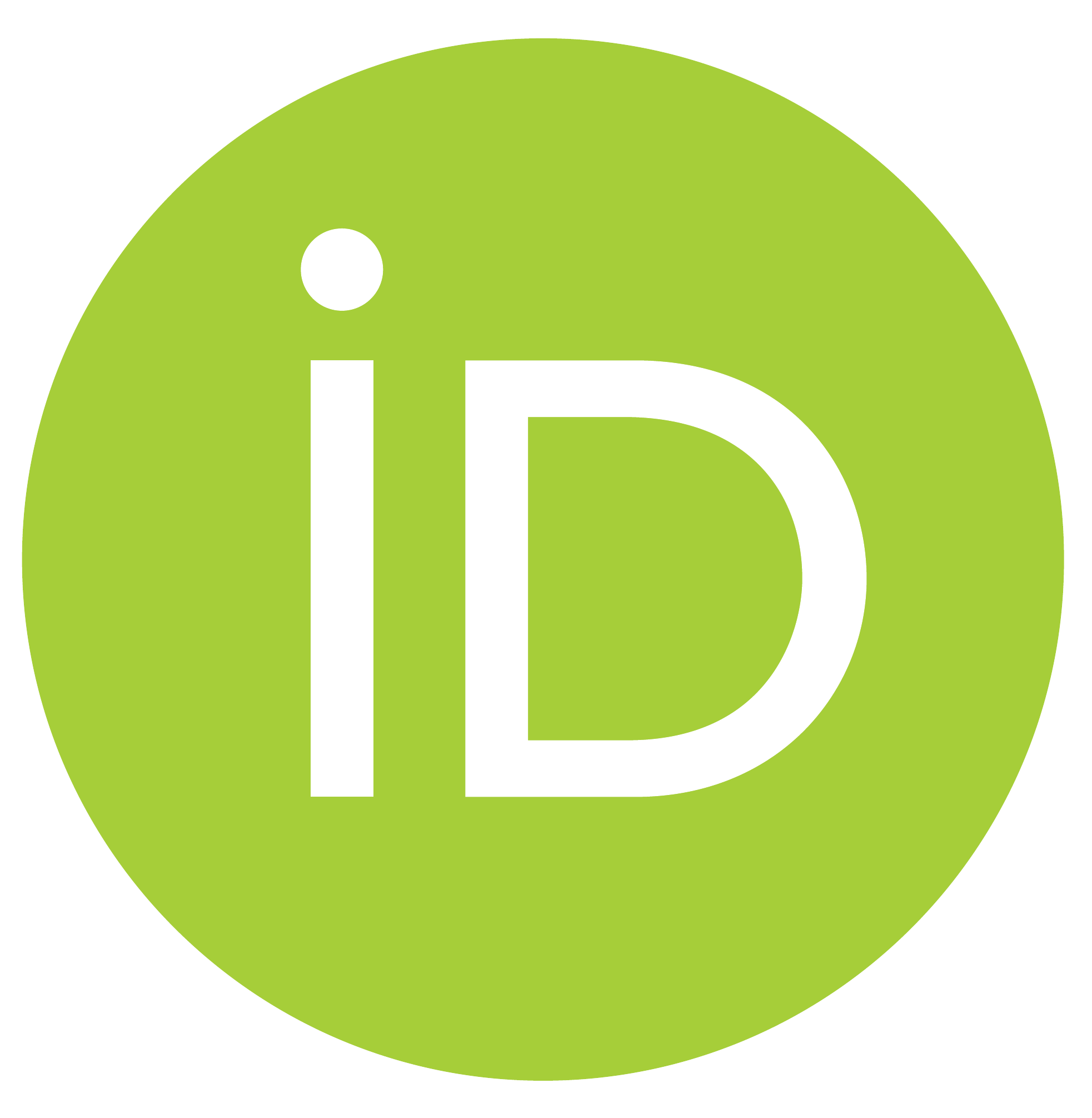}}}

\usepackage[titletoc,title]{appendix}
\usepackage{titlesec}
\titlespacing{\section}{0pt}{2.5ex}{1.5ex}
\titlespacing{\subsection}{0pt}{1.5ex}{1ex}
\titlespacing{\subsubsection}{0pt}{1.5ex}{1ex}
\titleformat{\section}{\large\bfseries\centering}{\thesection}{1em}{}
\titleformat{\subsection}[runin]{\bfseries}{\thesubsection.}{0.5em}{}[.\mbox{\ }]
\titleformat{\subsubsection}[runin]{\bfseries}{\thesubsubsection.}{0.4em}{}[.\mbox{\ }]

\usepackage{amsmath,amsfonts,amssymb,bm,mathrsfs,amsthm,amsxtra}
\numberwithin{equation}{section}
\newtheorem{theorem}{Theorem}[section]

\newtheorem{remark}{Remark}[section]

\usepackage{ulem}

\newcommand{\VERT}{\vert\kern-0.3ex\vert\kern-0.3ex\vert}

\def\p{\partial}

\usepackage{accents}
\makeatletter
\def\widebar{\accentset{{\cc@style\underline{\mskip10mu}}}}
\makeatother

\begin{document}
\title{\bf Well-posedness for moving interfaces in anisotropic plasmas\let\thefootnote\relax\footnotetext{
This research was carried out at the Sobolev Institute of Mathematics, under a state contract (project no. FWNF-2022-0008).
}
}

\author{
{\sc Yuri Trakhinin}\orcidicon{0000-0001-8827-2630}\thanks{Sobolev Institute of Mathematics, Koptyug av.~4, 630090 Novosibirsk, Russia. e-mail: trakhin@math.nsc.ru}
}

\date{\today}

\maketitle

\vspace{-6mm}

{\footnotesize
\noindent{\bf Abstract}:\quad
We study the local-in-time well-posedness for an interface that separates an anisotropic plasma from a vacuum. The plasma flow is governed by the ideal Chew--Goldberger--Low (CGL) equations, which are the simplest collisionless fluid model with anisotropic pressure. The vacuum magnetic and electric fields are supposed to satisfy the pre-Maxwell equations. The plasma and vacuum magnetic fields are tangential to the interface. This represents a nonlinear hyperbolic-elliptic coupled problem with a characteristic free boundary. By a suitable symmetrization of the linearized CGL equations we reduce the linearized free boundary problem to a problem analogous to that in isotropic magnetohydrodynamics (MHD). This enables us to prove the local existence and uniqueness of solutions to the nonlinear free boundary problem under the same non-collinearity condition for the plasma and vacuum magnetic fields on the initial interface required by Secchi and Trakhinin (Nonlinearity 27:105--169, 2014) in isotropic MHD.

\vspace{1mm}

\vspace{1mm}

\section{Introduction}\label{sec:intro}

Let $\Omega\subset \mathbb{R}^3$ be the reference domain occupied by a collisionless plasma and vacuum.
In the plasma region $\Omega^+(t)\subset \Omega$, the motion is governed by the following ideal Chew--Goldberger--Low (CGL) equations:
\begin{subequations}
	\label{MHD}
\begin{alignat}{2}
&    \p_t \rho +\nabla\cdot (\rho v )=0, \label{MHD1} \\
&\p_t (\rho v )+\nabla\cdot \left(\rho v \otimes v + (\tau -1) H \otimes H \right)+\nabla q =0, \label{MHD2} \\
&\p_t H-\nabla\times(v\times H)=0, \label{MHD3} \\
&\frac{d}{dt} \left(\frac{p_{\|}|H|^2}{\rho^3}\right) =0, \label{MHD4}\\
&\frac{d}{dt} \left(\frac{p_{\bot}}{\rho|H|}\right) =0, \label{MHD5}
\end{alignat}
\end{subequations}
 together with the divergence-free equation
\begin{align} \label{divH1}
  \nabla\cdot H=0.
\end{align}
Here density $\rho$, fluid velocity $v=(v_1,v_2,v_3)^{\mathsf{T}}$, magnetic field $H=(H_1,H_2,H_3)^{\mathsf{T}}$, parallel pressure $p_{\|}$, and perpendicular pressure $p_{\bot}$ are unknown functions of time $t$ and space variable $x=(x_1,x_2,x_3)$. We denote by
\[
\tau = \frac{p_{\|}-p_{\bot}}{|H|^2} \quad \mbox{and}\quad q=p_{\bot}+\tfrac{1}{2}|H|^2
\]
the anisotropy factor and the (perpendicular) total pressure respectively, and $d/dt= \partial_t +(v\cdot\nabla ) $ represents the material derivative. The system of equations \eqref{MHD} is closed for the primary unknown $U:=(\rho ,v,H,p_{\|},p_{\bot})^{\mathsf{T}}\in\mathbb{R}^9$ whereas \eqref{divH1} is the divergence constraint on the initial data for $U$.

The CGL equations named after Chew et al. \cite{CGL} are the simplest fluid model describing the motion of a collisionless plasma, i.e., plasma in which the mean free path for particle collisions is large compared to the Larmor radius. As we can see in \eqref{MHD2}, in the CGl model the pressure tensor has the anisotropic (gyrotropic) form ${\sf p} =p_{\bot}I +\tau H \otimes H$, where $I$ is the unit matrix of order 3. The CGL model is a gyrotropic approximation in the sense that in general the pressure tensor has also a non-gyrotropic part $\mathit{\Pi}$ connected with the finite Larmor radius (FLR) corrections: ${\sf p}=p_{\bot}I +\tau H \otimes H +\mathit{\Pi}$. The CGL model together with fluid models taking into account effects introduced by FLR corrections are discussed in full details in the nice survey \cite{Zank_et_al}. As was noted in \cite{Zank_et_al}, in recent years there has been an increased emphasis on pressure/temperature anisotropy effects, in particular, on the CGL model (see, e.g., \cite{Du,Melesh} and references therein) because the classical magnetohydrodynamic (MHD) fluid description does not satisfy all the needs of astrophysical applications requiring correct modelling of collisionless plasmas.

It is easy to see that in system \eqref{MHD} the classical CGL double adiabatic equations \eqref{MHD4} and \eqref{MHD5} can be equivalently replaced with the conservation laws
\begin{align}\label{MHD'}
\left\{
\begin{aligned}
&    \p_t (\rho s_{\|}) +\nabla\cdot (\rho s_{\|} v )=0,  \\
&    \p_t (\rho s_{\bot}) +\nabla\cdot (\rho s_{\bot} v )=0,
\end{aligned}
\right.
\end{align}
where
\[
s_{\|} =\frac{1}{3}\ln \left(\frac{p_{\|}|H|^2}{\rho^3}\right) \quad \mbox{and}\quad s_{\bot} = \frac{2}{3}\ln \left(\frac{p_{\bot}}{\rho|H|}\right)
\]
are so-called parallel and perpendicular entropies (see, e.g., \cite{Zank_et_al}). That is, equations \eqref{MHD1}--\eqref{MHD3}, \eqref{MHD'} form the system of 9 conservation laws. Mathematically, in \eqref{MHD'} instead of $s_{\|}$ and $s_{\bot}$ we could, of course, use any smooth functions of $p_{\|}|H|^2/\rho^3$ and $p_{\bot}/(\rho|H|)$ respectively.

Using \eqref{divH1} and the additional 10th conservation law (energy conservation)
\[
\p_t(\rho E + \tfrac{1}{2}|H|^2) +\nabla\cdot
\left( \rho Ev+p_{\bot}v +H\times(v\times H) +\tau (v\cdot H)H \right)=0
\]
which holds on smooth solutions of system \eqref{MHD} (see, e.g., \cite{Zank_et_al}), and following Godunov's symmetrization procedure \cite{God}, Blokhin and Krymskikh \cite{BK} have symmetrized the conservation laws \eqref{MHD1}--\eqref{MHD3}, \eqref{MHD'} in terms of a vector of canonical variables $Q=Q(U)$:
\[
A^0(Q )\partial_tQ+\sum_{j=1}^3A^j(Q )\partial_jQ=0.
\]
Here $E=\mathfrak{e}+\tfrac{1}{2} |v|^2$ is the total energy,
\[
\mathfrak{e}=\frac{p_{\bot}}{\rho}+\frac{p_{\|}}{2\rho}
\]
is the specific internal energy,
\[
Q = \left( \mathfrak{e} + \tfrac{p_{\|}}{\rho} -\tfrac{1}{2}|v|^2 -T_{\|}s_{\|}-T_{\bot}s_{\bot}\,, v\,, (1-\tau )H\,, T_{\|}\,,T_{\bot} \right)^{\mathsf{T}},
\]
$T_{\|}=p_{\|}/(\rho R)$ and $T_{\bot}=p_{\bot}/(\rho R)$ are the parallel and perpendicular temperatures, $R$ is the gas constant, and the symmetric matrices $A^{\alpha}(Q)$ ($\alpha =\overline{0,3}$) are written in \cite{BTmon}.

From  $Q$ we can return to the vector  $U$ of primary unknowns keeping the symmetry property (see \cite{BTmon}), i.e., we rewrite the CGL equations as the symmetric system
\begin{align}
	\label{MHD:vec}
	A_0^+(U)\p_t U+\sum_{j=1}^{3}A_j^+(U)\p_j U=0\qquad \textrm{in } \Omega^+(t),
\end{align}
The symmetric matrices $A_{\alpha}^+(U)$, which are rather cumbersome, are written in \cite{BTmon}. The symmetric system \eqref{MHD:vec} is hyperbolic if $A_0>0$. As was shown in \cite{BK} (see also \cite{BTmon}), the hyperbolicity condition $A_0^+>0$ holds provided that
\begin{equation}
\label{hypcond}
\rho >0,\quad
-1/a_p<\tau <1,
\end{equation}
where $a_p=T_{\bot}/T_{\|}=p_{\bot}/p_{\|}$ is the temperature anisotropy ratio \cite{Zank_et_al}. Moreover, we by default assume that $p_{\|}>0$, $p_{\bot}>0$, and the magnetic field $H\neq 0$.

In the vacuum region $\Omega^-(t)\subset\Omega$, for the vacuum magnetic field ${h}=(h_1,h_2,h_3)^{\mathsf{T}}$ and the vacuum electric field ${e}=(e_1,e_2,e_3)^{\mathsf{T}}$, we consider the pre-Maxwell equations
\begin{alignat}{3}
	& \nabla\times {h}=0,\qquad && \nabla \cdot {h}=0,
	\label{pre-M1a}
	\\
	&\nabla\times {e}=-\p_t h,\qquad && \nabla\cdot {e}=0,
	\label{pre-M1b}
\end{alignat}
where the displacement current is neglected from Maxwell's equations in vacuum as in non-relativistic CGL system. The vacuum electric field ${e}$ in \eqref{pre-M1a}--\eqref{pre-M1b} is a secondary variable, so that the dynamics in $\Omega^-(t)$ can be described by the elliptic (div-curl) system \eqref{pre-M1a}, or equivalently,
\begin{align}
	\label{pre-M:vec}
	\sum_{j=1}^{3}A_j^-\p_j {h}=0\qquad \textrm{in } \Omega^-(t),
\end{align}
where the constant matrices $A_1^-$, $A_2^-$, and $A_3^-$ are defined by
\begin{align}
	\nonumber  
	A_1^-:=\begin{pmatrix}
		0 & 0 &0 \\
		0 & 0 &-1 \\
		0 & 1 &0 \\
		1 & 0 &0
	\end{pmatrix},
	\quad
	A_2^-:=\begin{pmatrix}
		0 & 0 &1 \\
		0 & 0 &0 \\
		-1 & 0 &0 \\
		0 & 1 &0
	\end{pmatrix},
	\quad
	A_3^-:=\begin{pmatrix}
		0 & -1 &0 \\
		1 & 0 &0 \\
		0 & 0 &0 \\
		0 & 0 &1
	\end{pmatrix}.
\end{align}

As in \cite{ST14MR3151094}, for technical simplicity we assume that $\Omega^{\pm}(t)=\{x\in \Omega: x_1\gtrless \varphi(t,x') \}$ and the plasma--vacuum interface is given by the form of a graph
$$\Sigma(t):=\{x\in \Omega: x_1=\varphi(t,x')\}\quad \textrm{with \ } x'=(x_2,x_3),$$
where the interface function $\varphi$ is to be determined. As in \cite{ST14MR3151094}, we focus on the case of $\Omega=(-1,1)\times \mathbb{T}^2$ with boundaries $\Sigma^{\pm}:=\{\pm 1\}\times \mathbb{T}^2$,
where $\mathbb{T}^2$ denotes the $2$-torus and can be thought of as the unit square with periodic boundary conditions.
For the plasma--vacuum system the boundary conditions have the same form as in classical (isotropic) MHD \cite{ST14MR3151094}:
\begin{subequations}
	\label{BC1}
	\begin{alignat}{3}
		& q-\frac{1}{2}|{h}|^2=0,
		\quad \p_t \varphi=v\cdot N
		\qquad && \textrm{on } \Sigma(t),
		\label{BC1a}\\
		& H\cdot N=0,\quad
		{h}\cdot N=0
		&& \textrm{on } \Sigma(t),
		\label{BC1b}\\
		&H_1=0,\quad v_1=0 && \textrm{on } \Sigma^+,
		\label{BC1c}\\
		& {h}\times \mathbf{e}_1 =\bm{j}_{\rm c} &&\textrm{on } \Sigma^-,
		\label{BC1d}
	\end{alignat}
\end{subequations}
where $N:=(1,-\p_2\varphi,-\p_3\varphi)^{\mathsf{T}}$ is the normal to $\Sigma(t)$ and $\mathbf{e}_j:=(\delta_{1j},\delta_{2j},\delta_{3j})^{\mathsf{T}}$, $j=1,2,3$ with $\delta_{ij}$ being the Kronecker delta. The vector function $\bm{j}_{\rm c} $ represents a given surface current that forces oscillations onto the plasma--vacuum system. This model can be exploited for the analysis of waves in astrophysical plasmas, e.g., by mimicking the effects of excitation of MHD waves by an external plasma by means of a localized set of \lq\lq coils\rq\rq, when the response of the internal plasma is the main issue (e.g., in the problem of sunspot oscillations excited by sound waves in the photosphere; see \cite[\S 4.6]{GKP19} for a thorough discussion of the condition \eqref{BC1d}).

The second condition in \eqref{BC1a} means that the interface moves with the velocity of plasma particles. Conditions \eqref{BC1b} state that the plasma and vacuum magnetic fields are tangential to the interface. If on both sides of the interface we have plasmas, these conditions hold on a tangential discontinuity (current-vortex sheet). For the tangential discontinuity, the Rankine-Hugoniot jump conditions in anisotropic CGL plasmas \cite{BTmon,Hud} imply the jump condition $[q]=q^+|_{\Sigma}-g^-|_{\Sigma}=0$ for the perpendicular total pressure. That is, the first condition in \eqref{BC1a} appears as the limiting case of the jump condition $[q]=0$ when from one side of the discontinuity we have vacuum: $q^-=\frac{1}{2}|h|^2$. In other words, as in isotropic MHD, the first condition in \eqref{BC1a} comes from the balance of the normal stresses at the interface. At last, conditions \eqref{BC1c} are the standard perfectly conducting wall and impermeable conditions.

We supplement \eqref{MHD:vec} and \eqref{pre-M:vec}--\eqref{BC1} with the initial conditions
\begin{align} \label{IC1}
	\varphi|_{t=0}=\varphi_0,\qquad 	U|_{t=0}=U_0:=(\rho_0 ,v_0,H_0,p_{\|0},p_{\bot 0})^{\mathsf{T}},
\end{align}
where  $\|\varphi_0\|_{L^{\infty}(\mathbb{T}^2)}<1$.
Note that the vacuum magnetic field ${h}\in\mathbb{R}^3$ can be uniquely determined from the elliptic problem consisting of \eqref{pre-M:vec}, the second condition in \eqref{BC1b}, and \eqref{BC1d} when the interface function $\varphi$ is given.
It is worth mentioning that system \eqref{MHD:vec}, \eqref{pre-M:vec}--\eqref{IC1} is a nonlinear hyperbolic--elliptic coupled problem with a characteristic free boundary.

In \cite{T10MR2718711} two different well-posedness conditions we proposed for the linearized plasma--vacuum interface problem in classical ideal compressible MHD. The first one is the non-collinearity condition, stating that the magnetic fields on either side of the interface are not collinear:
\begin{align} \label{non-collinear}
|H\times h|\geq \kappa >0 \quad \textrm{on }\Sigma(t).
\end{align}
The second condition
\begin{align}\label{Taylor:q}
\bm{n}\cdot \nabla (q-\tfrac{1}{2}|h|^2)\leq -\kappa<0\quad \textrm{on }\Sigma(t)
\end{align}
is the MHD counterpart of the Rayleigh--Taylor sign condition, where $\bm{n}$ denotes the outward unit normal to the interface $\Sigma(t)$.

Based on the linear results in \cite{ST13MR3148595,T10MR2718711}, Secchi and Trakhinin \cite{ST14MR3151094} proved the first local well-posedness theorem for the (nonlinear) plasma--vacuum interface problem in classical MHD under condition \eqref{non-collinear} satisfied at the first moment. But, the proof of the local well-posedness of this problem under condition \eqref{Taylor:q} is still an open problem. At the same time, if the surface current $\bm{j}_{\rm c} \equiv 0$, the elliptic subproblem for $h$ has only zero solution $h\equiv 0$. For this case local well-posedness was showed by Trakhinin and Wang \cite{TW21MR4201624} under the Rayleigh--Taylor-type sign condition \eqref{Taylor:q}. In MHD, such case corresponds to the interface between a compressible liquid and vacuum because $h\equiv 0$ implies $q|_{\Sigma}=0$ (cf. \eqref{BC1a}). We note however that this is prohibited for the CGL model because $H\neq 0$ and $p_{\bot}>0$. At last, we refer the reader to \cite{HL14MR3187678,H17MR3614754,HL21MR4295166,GW19MR3980843,SWZ19MR3981394,MTT14MR3237563} and \cite{HL20MR4093862} respectively for well-posedness and ill-posedness results for plasma--vacuum interfaces in ideal incompressible MHD.

Our main goal in this paper is to extend the well-posedness result in \cite{ST14MR3151094} to the CGL model, i.e., to prove the local well-posedness of problem \eqref{MHD:vec}, \eqref{pre-M:vec}--\eqref{IC1} provided that the initial data satisfy the non-collinearity condition \eqref{non-collinear}. Fortunately, we do not need to repeat and adopt  to the CGL model all the long mathematical arguments in \cite{ST13MR3148595,ST14MR3151094}. Our main idea is finding a suitable symmetrization of the linearized CGL equations which enables us to reduce the linearized free boundary problem to a problem analogous to that in isotropic MHD. We also hope that this symmetrization could be useful for another boundary value problems in the CGL model, e.g., for current-vortex sheets or contact discontinuities \cite{BTmon}.

The plan of the rest of this paper is as follows. In Section \ref{sec:non}, we reduce the system \eqref{MHD:vec}, \eqref{pre-M:vec}--\eqref{IC1} to an equivalent fixed-boundary problem and state for it our main theorem on the local well-posedness under the non-collinearity condition \eqref{non-collinear}. In Section \ref{sec:lin}, we write down the linearized problem associated with the fixed-boundary problem from Section \ref{sec:non}. In Section \ref{sec:sym}, we obtain the announced symmetrization of the linearized CGL equations, and Section \ref{proof} is devoted to final remarks on the proof of our well-posedness theorem.

\section{Equivalent fixed-boundary problem and main result} \label{sec:non}

We reformulate the free boundary problem \eqref{MHD:vec}, \eqref{pre-M:vec}--\eqref{IC1} into an equivalent fixed-boundary problem by introducing
$U_{\sharp}(t,x):=U(t,\Phi(t,x),x')$ and $h_{\sharp}(t,x):=h(t,\Phi(t,x),x')$.
We choose the lifting function $\Phi$ as
\begin{align} \label{Phi:def}
	\Phi(t,x):=x_1+ \chi( x_1)\varphi(t,x'),
\end{align}
where $\chi\in C^{\infty}_0(-1,1)$ is the cut-off function that satisfies
$	\|\chi'\|_{L^{\infty}(\mathbb{R})} < 4/(\|\varphi_0\|_{L^{\infty}(\mathbb{T}^2)}+3)$ and equals to $1$ on a small neighborhood of the origin. See \cite{ST14MR3151094} for another change of variables, which can gain one half derivative for $\varphi$. Here, as in \cite{TW22c} and unlike \cite{ST14MR3151094}, we use the change of variables \eqref{Phi:def} for more technical simplicity (we refer to \cite{TW22c} for more details, in particular, for the restriction on $\chi'$, etc.).

After the change of variables \eqref{Phi:def} the free boundary problem \eqref{MHD:vec}, \eqref{pre-M:vec}--\eqref{IC1} is reduced to the following nonlinear fixed boundary problem:
\begin{subequations} \label{NP1}
	\begin{alignat}{2}
		\label{NP1a}
		&\mathbb{L}_+(U,\Phi) :=L_+(U,\Phi)U =0
		\quad & &\textrm{in  }   \Omega^+:=(0,1)\times \mathbb{T}^2,\\
		\label{NP1b}
		&\mathbb{L}_-({h},\Phi) :=L_-(\Phi){h} =0
		& &\textrm{in  }  		\Omega^-:=(-1,0)\times \mathbb{T}^2,\\
		\label{NP1c}
		&  \mathbb{B}(U,h,\varphi)=0
		&&\textrm{on  }     \Sigma^3\times \Sigma^+\times \Sigma^-,\\
		\label{NP1d}
		&U|_{t=0}=U_0,\qquad \varphi|_{t=0}=\varphi_0,  &  &
	\end{alignat}
\end{subequations}
where we have dropped the subscript $``\sharp"$ for convenience,
$\Sigma:=\{0\}\times \mathbb{T}^2$,
and
\begin{align}
	\label{L+:def}
	&L_+(U,\Phi):=A_0^+(U)\partial_t+\widetilde{A}_1^+(U,\Phi)\partial_1+ A_2^+(U)\partial_2
	+A_3^+(U)\partial_3,\\[0.5mm]
	\label{L-:def}
	&L_-(\Phi):=\widetilde{A}_1^-(\Phi)\partial_1+ A_2^-\partial_2+A_3^-\partial_3,\\[0.5mm]
\label{B:def}
&\mathbb{B}(U,h,\varphi):=
\begin{pmatrix}
 \p_t \varphi -v\cdot N\\  q-\tfrac{1}{2}|{h}|^2
 \\ {h}\cdot N\\  v_1\\ {h}\times \mathbf{e}_1-\bm{j}_{\rm c}
\end{pmatrix},
\end{align}
with $\widetilde{A}_1^-(\Phi):=  (A_1^--\partial_2\Phi A_2^--\partial_3\Phi A_3^-)/\partial_1\Phi$ and
\begin{align}
	\nonumber   
	&\widetilde{A}_1^+(U,\Phi):=
	\frac{1}{\partial_1\Phi}\big(A_1^+(U)-\partial_t\Phi A_0^+(U)-\partial_2\Phi A_2^+(U)-\partial_3\Phi A_3^+(U)\big).
\end{align}
In \eqref{NP1c}, we employ the notation $\Sigma^3\times \Sigma^+\times \Sigma^-$ to denote that the first three components of this vector equation are taken on $\Sigma$, the fourth one on $\Sigma^+$, and the fifth one on $\Sigma^-$. The equations for $H$ contained in \eqref{NP1a} can be written as
\begin{align} \label{H.bb:def}
\mathbb{H}(H,v,\Phi):=	
(\p_t^{\Phi}+v\cdot \nabla^{\Phi})H-(H\cdot \nabla^{\Phi})v+H\nabla^{\Phi}\cdot v=0\quad
\textrm{in }\Omega^+	,
\end{align}
where
\begin{align}
\nonumber
\p_t^{\Phi}:=\p_t-\frac{\p_t\Phi}{\p_1\Phi}\p_1,\
\nabla^{\Phi}:=(\p_1^{\Phi},\p_2^{\Phi},\p_3^{\Phi})^{\mathsf{T}},\
\p_1^{\Phi}:=\frac{\p_1}{\p_1\Phi},\
\p_j^{\Phi}:=\p_j-\frac{\p_j\Phi}{\p_1\Phi}\p_1
\end{align}
for $j=2,3$.
In the new variables, equation \eqref{divH1} and first conditions in \eqref{BC1b}--\eqref{BC1c} become
\begin{align}
	\label{inv1b}
\nabla^{\Phi}\cdot H =0\quad  \textrm{in  }    \Omega^+,\qquad
 H\cdot N=0\quad \textrm{on }    \Sigma, \qquad
 H_1=0\quad \textrm{on }    \Sigma^+,
\end{align}
which can be regarded as initial constraints,
meaning that they hold for $t>0$ as long as they are satisfied initially (see \cite{T09ARMAMR2481071} for the detailed proof).

For formulating of the main result of this paper, which is a local existence and uniqueness theorem for problem \eqref{NP1}, we need to introduce the anisotropic weighted Sobolev spaces \cite{C07MR2289911,S00MR1779863}. We denote
\begin{align}
\nonumber
\mathrm{D}_*^{\alpha}:=\p_t^{\alpha_0} (\sigma \p_1)^{\alpha_1}\p_2^{\alpha_2}  \p_3^{\alpha_3} \p_1^{\alpha_{4}}\qquad\textrm{for }  \alpha:=(\alpha_0,\ldots,\alpha_{4})\in\mathbb{N}^{5},	
\end{align}
where $\sigma=\sigma(x_1)$ is a positive $C^{\infty}$--function on $(0,1)$ such that $\sigma(x_1)=x_1$ in a  neighbourhood of the origin and $\sigma(x_1)=1-x_1$ in a  neighbourhood of $x_1=1$.
For $m\in\mathbb{N}$ and $I\subset \mathbb{R}$, the anisotropic Sobolev space $H_*^{m}(I\times \Omega^+)$ is defined as
\begin{align*}
	&H_*^m(I\times\Omega^+):=
	\{ u\in L^2(I\times\Omega^+):\, \mathrm{D}_*^{\alpha} u\in L^2(I\times\Omega^+)
	\textrm{ for } \langle \alpha \rangle\leq m  \},
\end{align*}
and equipped with the norm ${\|}\cdot{\|}_{H^m_*(I\times \Omega^+)}$, where
\begin{align}  \nonumber 
	\langle \alpha \rangle :=\sum_{i=0}^3\alpha_i +2\alpha_{4},\qquad	
	{\|}u{\|}_{H^m_*(I\times \Omega^+)}^2:=
	\sum_{\langle \alpha\rangle\leq m} \|\mathrm{D}_*^{\alpha} u\|_{L^2(I \times\Omega^+)}^2.
\end{align}
By definition, $H^m(I\times \Omega^+)\hookrightarrow H_*^m(I\times \Omega^+) \hookrightarrow H^{\lfloor m/2\rfloor}(I\times \Omega^+) $ for all  $m\in \mathbb{N}$ and  $I\subset \mathbb{R}$,
where $\lfloor s \rfloor$ denotes the floor function of $s\in\mathbb{R}$ that maps $s$ to the greatest integer less than or equal to $s$.
We refer to \cite{MST09MR2604255,S00MR1779863,C07MR2289911} and references therein for an extensive study of anisotropic Sobolev spaces.

To present our main theorem, we also need to introduce the compatibility conditions on the initial data. However, the process of introduction of the compatibility conditions for problem \eqref{NP1} totally coincides with that described in \cite{ST14MR3151094,TW22c} for the counterpart of \eqref{NP1} in isotropic MHD.
This is why we can just refer to \cite{ST14MR3151094,TW22c}. For the reader's convenience, we only note here that the initial vacuum magnetic field $h_0 $ is not given independently from the other initial data because it is uniquely determined by the div-curl system
\begin{align}
	\nonumber
	L_-(\Phi_0)h_0 =0\ \ \textrm{in  }  \Omega^-,\quad  \ \
	h_0\cdot N_0=0 \ \ \textrm{on  }    \Sigma,\quad  \ \
	h_0\times \mathbf{e}_1=\bm{j}_{\rm c}(0)\ \ \textrm{on  }    \Sigma^-,
\end{align}
where $L_-$ is the operator given by \eqref{L-:def} and $N_0:=(1,-\p_2\varphi_0,-\p_3\varphi_0)^{\mathsf{T}}$.

We are now in a position to state the main result of this paper.

\begin{theorem}
	\label{thm1}
Assume that  $\bm{j}_{\rm c}\in H^{m+3/2}([0,T_0]\times \Sigma^-)$ for some $T_0>0$ and $m\in\mathbb{N}$ with $m\geq 20$.
Assume further that the initial data $(U_0,\varphi_0)\in H^{m+3/2}(\Omega^+)\times H^{m+2}(\mathbb{T}^{2})$
satisfy $\|\varphi_0\|_{L^{\infty}(\mathbb{T}^2)}<1$, the constraints \eqref{inv1b}, the hyperbolicity conditions \eqref{hypcond},
\begin{align}\label{hyp}
\rho_0 \geq \delta_1>0,\quad \frac{p_{\bot 0}-p_{\|0}}{|H_0|^2} +1 \geq \delta_2>0,\quad p_{\|0} +\frac{p_{\bot 0}(p_{\|0}-p_{\bot 0})}{|H_0|^2} \geq \delta_3>0,
\end{align}
the condition
\begin{align}\label{hyp-lin}
6p_{\| 0}- p_{\bot 0}\geq \delta_4>0,
\end{align}
the default requirements  $p_{\| 0}\geq \delta_5>0$, $p_{\bot 0}\geq \delta_6>0$, $|H_0|\geq \delta_7>0$, the compatibility conditions up to order $m$, and the non-collinearity condition
\begin{align}
\label{non-col}
|H_0\times h_0|\big|_{\Sigma}\geq \delta_0>0
\end{align}
for some fixed constants $\delta_k$ ($k=\overline{0,4}$). Then problem \eqref{NP1} admits a unique solution $(U,h, \varphi)$ in $H_*^{m-9}([0,T]\times\Omega^+)\times H^{m-9}([0,T]\times\Omega^-)\times H^{m-9}([0,T]\times\mathbb{T}^{2})$ for some $0<T\leq T_0$.
\end{theorem}

\begin{remark}
According to the local existence results in \cite{K75MR0407477,V-H} for general symmetric hyperbolic systems, the Cauchy problem for the CGL equations \eqref{MHD} allows smooth solutions within a short time if the initial data satisfy the hyperbolicity conditions \eqref{hypcond}/\eqref{hyp}. Clearly, conditions \eqref{hypcond} satisfied for a background constant solution prevent the ill-posedness of the Cauchy problem for the corresponding constant coefficient linearized CGL equations, in particular, the so-called firehose  and mirror instabilities \cite{Zank_et_al}  taking place if $\tau >1$ and $a_p >6 (1 +1/\beta_{\bot})$ respectively, where $\beta_{\bot} = 2p_{\bot}/|H|^2$ is the perpendicular plasma beta (one can easily check that the inequality $\tau> -1/a_p$ in \eqref{hypcond} prevents mirror instability).
\label{r1}
\end{remark}

\begin{remark}
In view of the inequality $\tau> -1/a_p$ in \eqref{hypcond}, our restriction $a_p<6$ on the initial data in \eqref{hyp-lin} is automatically satisfied for $\beta_{\bot}>2/5$. Otherwise, it is indeed an additional requirement on the initial data. We believe that it could be neglected if we were able to get a priori estimates for the linearization of problem \eqref{NP1} without using the symmetrization of the linearized CGL equations in Section \ref{sec:sym} alternative to the linearization of the symmetric system \eqref{MHD:vec} with really cumbersome matrices. On the other hand, even if condition \eqref{hyp-lin} can be dropped, we hope that the symmetrization proposed in Section \ref{sec:sym} could be useful for applying the energy method to another boundary value problems for the CGL system.
\label{r2}
\end{remark}

\begin{remark}
Since in this paper we use for technical simplicity the change of variables \eqref{Phi:def}, which does not gain one half derivative for $\varphi$ compared to the change in \cite{ST14MR3151094}, Theorem \ref{thm1} is formulated in the form similar to the theorem proved in \cite{TW22c} for the plasma--vacuum interface problem in isotropic MHD with non-zero surface tension.
\label{r3}
\end{remark}


\section{Linearized problem} \label{sec:lin}

The proof of the existence and uniqueness of solutions to a nonlinear problem is often relied on the analysis of the linearized problem. Moreover, if for the linearized problem we can only deduce a priori estimates with a loss of derivatives from the source terms and coefficients, then the classical fixed-point argument cannot be applied. This difficulty can be sometimes overcome by using the Nash--Moser method (see, e.g., \cite{S16MR3524197} and references therein). To this end, one needs to perform a ``genuine'' linearization when we keep all the lower-order terms which are usually dropped while applying the fixed-point argument. Following arguments in \cite{ST13MR3148595,ST14MR3151094,TW22c} for the plasma-vacuum interface problem in classical MHD, we below describe such a linearization of problem \eqref{NP1}.

Consider the basic state $(\mathring{U}(t,x),\mathring{{h}}(t,x),\mathring{\varphi}(t,x'))$ which is
is a sufficiently smooth vector-function defined on $\Omega_{T}^+\times \Omega_{T}^-\times\Sigma_{T}$, where $\Omega_T^{\pm}:=(-\infty,T)\times \Omega^{\pm}$,
and $\Sigma_T:=(-\infty,T)\times \Sigma$. The assumptions on the basic state are totally analogous to those in \cite{ST14MR3151094,TW22c}. In particular, we assume that it satisfies the hyperbolicity conditions ({\it cf.}~\eqref{hyp})
\begin{align}\label{hypl}
\mathring{\rho} \geq \frac{\delta_1}{2}>0,\quad \frac{\mathring{p}_{\bot }-\mathring{p}_{\|}}{|\mathring{H}|^2} +1 \geq \frac{\delta_2}{2}>0,\quad p_{\|} +\frac{\mathring{p}_{\bot }(\mathring{p}_{\|}-\mathring{p}_{\bot })}{|\mathring{H}|^2} \geq \frac{\delta_3}{2}>0
\end{align}
(together with $\mathring{p}_{\| }\geq \delta_5/2>0$, $\mathring{p}_{\bot }\geq \delta_6/2>0$, $|\mathring{H}|\geq \delta_7/2>0$), the condition ({\it cf.}~\eqref{hyp-lin})
\begin{align}\label{hyp-lin'}
6- \frac{\mathring{p}_{\bot }}{\mathring{p}_{\| }}\geq \frac{\delta_4}{2}>0,
\end{align}
some regularity assumptions as in \cite{ST14MR3151094,TW22c} and
\begin{alignat}{3}
\label{bas1g}
&\mathbb{H}(\mathring{H},\mathring{v},\mathring{\Phi})=0	
	\qquad &&\textrm{in }\Omega^+_T,\\	
	\label{bas1d}
	&\p_t \mathring{\varphi}=\mathring{v}\cdot \mathring{N},
	\quad \mathring{{h}}\cdot \mathring{N}=0
	\qquad &&\textrm{on }\Sigma_T,\\	
	\label{bas1e}
	&  \p_1 \mathring{{h}}\cdot \mathring{N}+\p_2\mathring{h}_2+\p_3\mathring{h}_3=0
	\quad &&\textrm{on }\Sigma_T,
	\\
	\label{bas1f}
	&
	\mathring{v}_1=0
	\quad  \textrm{on }\Sigma_T^+,
	\qquad
	\mathring{{h}}\times \mathbf{e}_1 =\bm{j}_{\rm c}
	\quad &&\textrm{on }\Sigma_T^-,
\end{alignat}
where $\mathbb{H}$ is the operator defined in \eqref{H.bb:def},  $\mathring{U}=(\mathring{\rho} ,\mathring{v},\mathring{H},\mathring{p}_{\|},\mathring{p}_{\bot})^{\mathsf{T}}\in\mathbb{R}^9$,
$\mathring{{h}}=(\mathring{h}_1,\mathring{h}_2,\mathring{h}_3)^{\mathsf{T}}\in\mathbb{R}^3,$ $\mathring{{v}}=(\mathring{v}_1,\mathring{v}_2,\mathring{v}_3)^{\mathsf{T}}\in\mathbb{R}^3,$
$\mathring{\Phi}(t,x):=x_1+\mathring{\Psi}(t,x)$ with
$\mathring{\Psi}(t,x):=\chi(x_1)\mathring{\varphi}(t,x'),$
$\mathring{N}:=(1,-\p_2\mathring{\Phi},-\p_3\mathring{\Phi})^{\mathsf{T}}$, and $\Sigma_T^{\pm}:=(-\infty,T)\times \Sigma^{\pm}$.
It follows from \eqref{bas1g} that the identities
\begin{align} \label{bas1h}
\nabla^{\mathring{\Phi}}\cdot \mathring{H}\big|_{\Omega^+_T} =0\qquad
\mathring{H}\cdot \mathring{N}\big|_{\Sigma_T}=0,\qquad
\mathring{H}_1\big|_{\Sigma_T^+}=0
\end{align}
are satisfied if they hold at the initial time (see \cite{T09ARMAMR2481071} for the proof). As such, we require that the conditions \eqref{bas1h} are satisfied at $t=0$.
Moreover, we assume that the non-collinearity condition holds for the basic state ({\it cf.}~\eqref{non-col}):
\begin{align}
	\label{bas:non-col}
\big|\mathring{H}\times\mathring{h}\big|\geq \frac{\delta_0}{2}>0	\qquad \textrm{on }\Sigma_{T}.
\end{align}

Introduce the good unknowns of Alinhac \cite{A89MR976971}:
\begin{align} \label{good}
	\dot{U}:=U-\frac{\Psi}{\partial_1 \mathring{\Phi}}\partial_1\mathring{U},
	\quad \dot{{h}}:={h}-\frac{\Psi}{\p_1 \mathring{\Phi}}\p_1 \mathring{{h}},
\end{align}
where $\Psi(t,x):=\chi(x_1)\psi(t,x')$.
Then the linearized operators for equations \eqref{NP1a}--\eqref{NP1b} around the basic state $(\mathring{U},\mathring{{h}},\mathring{\varphi})$ read
\begin{align}
	\mathbb{L}'_+\big(\mathring{U},\mathring{\Phi}\big)(U,\Psi)
	:=&\left.\frac{\mathrm{d}}{\mathrm{d}\theta}
	\mathbb{L}_{+}\big(\mathring{U} +{\theta}U ,
	\mathring{\Phi} +{\theta}\Psi \big)\right|_{\theta=0}
	\nonumber%
	\\
	=& \,L_{+}\big(\mathring{U},\mathring{\Phi}\big)\dot{U}+\mathcal{C}_{+}(\mathring{U},\mathring{\Phi})\dot{U}  +\frac{\Psi}{\partial_1\mathring{\Phi}}
	\partial_1\mathbb{L}_+(\mathring{U},\mathring{\Phi} ),
	\label{Alinhac1}
	\\[0.5mm]
	\mathbb{L}_{-}'\big(\mathring{{h}},\mathring{\Phi}\big)({h},\Psi)
	:=& \left.\frac{\mathrm{d}}{\mathrm{d}\theta} \mathbb{L}_{-}\big(\mathring{{h}} +{\theta}{h} , \mathring{\Phi} +{\theta}\Psi \big)\right|_{\theta=0}
	=L_-\big(\mathring{\Phi}\big)\dot{{h}} +\frac{\Psi}{\partial_1\mathring{\Phi}}\partial_1\mathbb{L}_-(\mathring{{h}},\mathring{\Phi} ),
	\label{Alinhac2}
\end{align}
where $L_{\pm}$ are the operators defined in \eqref{L+:def}--\eqref{L-:def} and
\begin{align} \nonumber 
	\mathcal{C}_+({U},{\Phi})U:=
	\sum_{k=1}^{8}V_k\Bigg(
	\frac{\p \widetilde{A}^{+}_1}{\p {U_k}}({U},{\Phi}) \partial_1 {U}
	+\sum_{i=0,2,3}\frac{\p A^{+}_i}{\p {U_k}}({U}) \partial_i {U}
	\Bigg).
\end{align}
For the boundary operator $\mathbb{B}$ defined by \eqref{B:def}, we have
\begin{align}
	\nonumber
&\mathbb{B}'(\mathring{U} ,\mathring{h}, \mathring{\varphi})(U,h,\psi)
:=  \left.\frac{\mathrm{d}}{\mathrm{d}\theta}
	\mathbb{B}\big(\mathring{U} +\theta U ,\,\mathring{h} +\theta h ,\,\mathring{\varphi} +\theta \psi \big) \right|_{\theta=0} \\[1mm]
&\ \ =	\left((\p_t + \mathring{v}' \cdot \mathrm{D}_{x'}) \psi-v\cdot\mathring{N},
\;
	p_{\bot}+\mathring{H}\cdot H-\mathring{h}\cdot h,\;
	h\cdot\mathring{N}-  \mathring{h}' \cdot \mathrm{D}_{x'} \psi,\;
	v_1,\;
	h\times \mathbf{e}_1		
\right)^{\mathsf{T}},
\nonumber 
\end{align}
where we denote $z':=(z_2,z_3)^{\mathsf{T}}$ for any vector $z:=(z_1,z_2,z_3)^{\mathsf{T}}$.

Dropping the last terms in \eqref{Alinhac1}--\eqref{Alinhac2}, we get the following effective linear problem for the good unknowns \eqref{good}:
\begin{subequations} \label{ELP1}
	\begin{alignat}{3}
		&\mathbb{L}'_{e+}(\mathring{U}, \mathring{\Phi}) \dot{U}
		:=L_+(\mathring{U}, \mathring{\Phi})\dot{U}
		+\mathcal{C}_+( \mathring{U},\mathring{\Phi})\dot{U} =f^+
		&\quad  &\textrm{in } \Omega^+_T,
		\label{ELP1a}\\
		&L_-( \mathring{\Phi}) \dot{{h}}=f^-
		&\quad  &\textrm{in } \Omega^-_T,
		\label{ELP1b}\\
		\label{ELP1c}   &
		\mathbb{B}'_{e}(\mathring{U}, \mathring{{h}}, \mathring{\varphi}) (\dot{U},\dot{{h}},\psi)
		=g
		&&\textrm{on } \Sigma_T^3\times \Sigma_T^+\times \Sigma_{T}^-,\\
		&(\dot{U},\psi)\big|_{t<0}=0,\qquad \dot{{h}}\big|_{t<0}=0,
		&&
		\label{ELP1d}
	\end{alignat}
\end{subequations}
where
\begin{align}
\label{B'e:def}	
\mathbb{B}'_{e}(\mathring{U}, \mathring{{h}}, \mathring{\varphi}) (\dot{U},\dot{{h}},\psi)
:=\begin{pmatrix}
	(\p_t + \mathring{v}' \cdot \mathrm{D}_{x'}+\mathring{b}_1) \psi-\dot{v}\cdot\mathring{N}\\[0.5mm]
\dot{p}_{\bot}+\mathring{H}\cdot \dot{H}-\mathring{h}\cdot \dot{h}+\mathring{b}_2  \psi\\[0.5mm]
\dot{{h}}\cdot \mathring{N}-\mathrm{D}_{x'}\cdot(\mathring{h}'  \psi)\\[0.5mm]
\dot{v}_1\\[0.5mm]		
\dot{h}\times \mathbf{e}_1		
\end{pmatrix}
\end{align}
with $\mathring{b}_1:=-\p_1 \mathring{v}\cdot \mathring{N}$ and
$\mathring{b}_2:=\p_1\mathring{p}_{\bot}+\mathring{H}\cdot\p_1 \mathring{H}-\mathring{{h}}\cdot \p_1\mathring{{h}}$, definitions \eqref{good}, and constraint \eqref{bas1e}. The last terms in \eqref{Alinhac1}--\eqref{Alinhac2} should be considered as error terms at each Nash--Moser iteration step (see \cite{ST14MR3151094,TW22c} for more details). The source terms $f^{\pm}$ and $g$ are supposed to vanish in the past.  We consider the case of zero initial data because the nonlinear problem can be reduced to it by construction of a so-called approximate solution (see, e.g., \cite{ST14MR3151094,TW22c}).

\section{Symmetrization of the linearized CGL equations} \label{sec:sym}

To apply the energy method to the linear problem \eqref{ELP1}, we need to have symmetric matrices in the operator $L_+$ described in \eqref{L+:def}. At the same time, the symmetric matrices in system \eqref{MHD:vec} found in \cite{BK}  are so cumbersome (see also \cite{BTmon}) that it makes them difficult to be used for deriving a priori estimates for problem \eqref{ELP1}. To avoid this difficulty we propose an elementary (algebraic) symmetrization of the linearized CGL equations, which does not rely on the Godunov's symmetrization procedure.

We forget for a moment about our initial-boundary value problem  \eqref{ELP1} and just consider the linearization of the CGL equation in the whole space $\mathbb{R}^3$ about a basic state $\mathring{U}=(\mathring{\rho} ,\mathring{v},\mathring{H},\mathring{p}_{\|},\mathring{p}_{\bot})^{\mathsf{T}}$.  We first rewrite system \eqref{MHD} in the following quasilinear form:
\begin{equation}
\p_t U + \sum_{j=1}^{3}B_j(U)\p_j U=0,
\label{quasi}
\end{equation}
where the matrices $B_j(U)$ are not symmetric and can be written down if necessary. Note that the divergence constraint \eqref{divH1} is used while writing down system \eqref{quasi}. Then the linearization of \eqref{quasi} about $\mathring{U}$ reads:
\begin{align}\label{linMHD}
\left\{
\begin{aligned}
&    D(\mathring{v})\rho+\mathring{\rho}\,\nabla\cdot v + {\rm z.o.t.}=0,  \\
&    \mathring{\rho}\,D(\mathring{v})v +\mathring{b} \big(\mathring{b}\cdot \big\{\nabla (p_{\|}-p_{\bot})-2\mathring{\tau}(\mathring{H}\cdot\nabla )H\big\}\big) + (\mathring{\tau}-1)(\mathring{H}\cdot\nabla )H +\nabla q + {\rm z.o.t.}=0,\\
& D(\mathring{v})H- (\mathring{H}\cdot\nabla )v + \mathring{H}\,\nabla\cdot v + {\rm z.o.t.}=0,\\
& D(\mathring{v})p_{\|} +\mathring{p}_{\|}\,\nabla\cdot v +2\mathring{p}_{\|}\big(\mathring{b}\cdot(\mathring{b}\cdot\nabla )v\big) + {\rm z.o.t.}=0,\\
& D(\mathring{v})p_{\bot} +2\mathring{p}_{\bot}\,\nabla\cdot v -\mathring{p}_{\bot}\big(\mathring{b}\cdot(\mathring{b}\cdot\nabla )v\big) + {\rm z.o.t.}=0,
\end{aligned}
\right.
\end{align}
where $U=(\rho ,v,H,p_{\|},p_{\bot})^{\mathsf{T}}$ is now the vector of perturbations,
\[
D(\mathring{v}) =\partial_t +(\mathring{v}\cdot\nabla ),\quad \mathring{b}=\frac{\mathring{H}}{|\mathring{H}|},\quad \mathring{\tau} = \frac{\mathring{p}_{\|}-\mathring{p}_{\bot}}{|\mathring{H}|^2}, \quad  q=p_{\bot}+\mathring{H}\cdot H,
\]
and z.o.t. are zero-order (nondifferential) terms which are of no interest for our current goals.

Let us introduce the new unknown
\[
P:= \tfrac{1}{2}p_{\bot}-p_{\|}+\mathring{\tau}(\mathring{H}\cdot H).
\]
The introduction of this unknown is prompted by some functions used in \cite{BT93} for studying by the energy method the 2D stability of rectilinear shock waves in the CGL model for the special cases when the background constant magnetic field is parallel or perpendicular to the shock front. Let the new vector-valued unknown function be
\begin{equation}
\label{change}
V= (p_{\bot},v,H,P,s_{\|})^{\mathsf{T}} = J(\mathring{U})U,
\end{equation}
where $s_{\|}$ is now the perturbation of the parallel entropy:
\[
s_{\|}:=\frac{p_{\|}}{3\mathring{p}_{\|}}-\frac{\rho}{\mathring{\rho}}+\frac{2(\mathring{H}\cdot H)}{3|\mathring{H}|^2}
\]
(instead of the $s_{\|}$ we could also take the perturbation of the perpendicular entropy), and the nonsingular matrix $J(\mathring{U})$ can be easily written down. We then rewrite \eqref{linMHD} in terms of the new unknown $V$ as follows:
\begin{align}\label{linMHD'}
\left\{
\begin{aligned}
& \frac{1}{2\mathring{p}_{\bot}}D(\mathring{v})p_{\bot} +\nabla\cdot v -\tfrac{1}{2}\big(\mathring{b}\cdot(\mathring{b}\cdot\nabla )v\big) + {\rm z.o.t.}=0,\\[3pt]
& \mathring{\rho}\,D(\mathring{v})v -\mathring{b} \big(\mathring{b}\cdot \big\{\nabla \big(\tfrac{1}{2}p_{\bot} +P\big)+\mathring{\tau}(\mathring{H}\cdot\nabla )H\big\}\big) + (\mathring{\tau}-1)(\mathring{H}\cdot\nabla )H +\nabla q + {\rm z.o.t.}=0,\\[3pt]
& (1-\mathring{\tau})D(\mathring{v})H +\mathring{\tau} \mathring{b}\big(\mathring{b}\cdot D(\mathring{v})H\big) +(\mathring{\tau}-1) (\mathring{H}\cdot\nabla )v \\[3pt]
& \qquad\qquad\qquad\qquad\qquad
-\mathring{\tau} \mathring{b}\big(\mathring{b}\cdot (\mathring{H}\cdot\nabla )v \big)
+ \mathring{H}\,\nabla\cdot v + {\rm z.o.t.}=0,\\[3pt]
& \frac{2}{6\mathring{p}_{\|}-\mathring{p}_{\bot}}\,D(\mathring{v})P  -\mathring{b}\cdot(\mathring{b}\cdot\nabla )v  + {\rm z.o.t.}=0,\\[3pt]
& D(\mathring{v})s_{\|} + {\rm z.o.t.}=0.
\end{aligned}
\right.
\end{align}
For obtaining equations \eqref{linMHD'} we have, in particular, left multiplied the equation for $H$ in \eqref{linMHD} by the symmetric matrix
$$
\mathring{\mathcal{B}}=(1-\mathring{\tau})I+\mathring{\tau}\mathring{b}\otimes\mathring{b}.
$$

Equations \eqref{linMHD'} form the symmetric hyperbolic system
\begin{equation}\label{linsym}
\mathcal{A}_0(\mathring{U})\p_t V+\sum_{j=1}^{3}\mathcal{A}_j(\mathring{U})\p_j V +\mathcal{A}_4(\mathring{U})V=0,
\end{equation}
where the block diagonal matrix
$\mathcal{A}_0(\mathring{U})={\rm diag}\,(1/(2\mathring{p}_{\bot}),\mathring{\rho}I,\mathring{\mathcal{B}},2/(6\mathring{p}_{\|}-\mathring{p}_{\bot}),1)$
is positive definite thanks to the assumptions $\mathring{\tau}<1$ (cf. \eqref{hypl}) and \eqref{hyp-lin'},
\[
\mathcal{A}_j(\mathring{U})=
\begin{pmatrix}
\dfrac{\mathring{v}_j}{2\mathring{p}_{\bot}} &\mathbf{e}_j^{\mathsf{T}}-\tfrac{1}{2}\mathring{b}_j\mathring{b}^{\mathsf{T}}\; & 0 & 0\; & 0 \\[9pt]
\mathbf{e}_j-\tfrac{1}{2}\mathring{b}_j\mathring{b}\; & \mathring{\rho}\mathring{v}_jI & \mathbf{e}_j\otimes \mathring{H} -\mathring{H}_j\mathring{\mathcal{B}}\quad &  -\mathring{b}_j\mathring{b}\;& 0 \\[6pt]
0 & \mathring{H} \otimes\mathbf{e}_j  -\mathring{H}_j\mathring{\mathcal{B}}\quad & \mathring{v}_j\mathring{\mathcal{B}} \;& 0 & 0 \\[9pt]
0 & -\mathring{b}_j\mathring{b}^{\mathsf{T}} & 0 & \dfrac{2\mathring{v}_j}{6\mathring{p}_{\|}-\mathring{p}_{\bot}}\;\, & 0 \\[9pt]
0 & 0 & 0 & 0 & \mathring{v}_j
\end{pmatrix}
,
\]
and the concrete form of the matrix $\mathcal{A}_4(\mathring{U})$ is of no interest.

\section{Proof of Theorem \ref{thm1}} \label{proof}

In view of \eqref{change}, \eqref{linsym}, for deriving a priori estimates for the linear problem \eqref{ELP1} we can now use the system
\begin{equation}\label{linsym'}
\mathcal{A}_0(\mathring{U})\p_t \dot{V}+\widetilde{\mathcal{A}}_1(\mathring{U},\mathring{\Phi})\p_1 \dot{V} +\mathcal{A}_2(\mathring{U})\p_2 \dot{V}+\mathcal{A}_3(\mathring{U})\p_3 \dot{V}+\mathcal{A}_4(\mathring{U})\dot{V}=F^+\quad \textrm{in } \Omega^+_T
\end{equation}
following from \eqref{ELP1a}, where
\[
\dot{V}:=J(\mathring{U})\dot{U},\qquad
\widetilde{\mathcal{A}}_1^+(\mathring{U},\mathring{\Phi}):=
	\frac{1}{\partial_1\mathring{\Phi}}\big(\mathcal{A}_1^+(\mathring{U})-\partial_t\mathring{\Phi} \mathcal{A}_0^+(\mathring{U})-\partial_2\mathring{\Phi} \mathcal{A}_2^+(\mathring{U})-\partial_3\mathring{\Phi} \mathcal{A}_3^+(\mathring{U})\big),
\]
and the new source term $F^+= \mathcal{A}_0(\mathring{U})\big(J(\mathring{U})\big)^{-1}\big(A_0(\mathring{U})\big)^{-1}f^+$. The crucial role is then played by the boundary matrix $\widetilde{\mathcal{A}}_1(\mathring{U},\mathring{\Phi})$ calculated on the boundary $\Sigma_T$:
\[
\widetilde{\mathcal{A}}_1(\mathring{U},\mathring{\Phi})=
\begin{pmatrix}
\dfrac{\mathring{m}}{2\mathring{p}_{\bot}} &\mathring{N}^{\mathsf{T}}-\tfrac{1}{2}(\mathring{b}\cdot \mathring{N})\mathring{b}^{\mathsf{T}}\; & 0 & 0\; & 0 \\[9pt]
\mathring{N}-\tfrac{1}{2}(\mathring{b}\cdot \mathring{N})\mathring{b}\; & \mathring{\rho}\mathring{m}I & \mathring{N}\otimes \mathring{H} -(\mathring{H}\cdot \mathring{N})\mathring{\mathcal{B}}\quad &  -(\mathring{b}\cdot \mathring{N})\mathring{b}\;& 0 \\[6pt]
0 & \mathring{H} \otimes\mathring{N}-(\mathring{H}\cdot \mathring{N})\mathring{\mathcal{B}}\quad & \mathring{m}\mathring{\mathcal{B}} \;& 0 & 0 \\[9pt]
0 & -(\mathring{b}\cdot \mathring{N})\mathring{b}^{\mathsf{T}} & 0 & \dfrac{2\mathring{m}}{6\mathring{p}_{\|}-\mathring{p}_{\bot}}\;\, & 0 \\[9pt]
0 & 0 & 0 & 0 & \mathring{m}
\end{pmatrix}
\quad \textrm{on } \Sigma_T,
\]
where $\mathring{m}:=\mathring{v}\cdot \mathring{N}-\p_t\mathring{\varphi}$. It follows from \eqref{bas1d} that
\[
\widetilde{\mathcal{A}}_1(\mathring{U},\mathring{\Phi})=
\begin{pmatrix}
0 &\mathring{N}^{\mathsf{T}} & 0 & 0 & \;0 \\
\mathring{N}\; & 0 & \mathring{N}\otimes \mathring{H}\;  &  0& \;0 \\
0 & \mathring{H} \otimes\mathring{N} \;&  0 & 0 &  \;0\\
0 & 0 & 0 & 0 &\; 0 \\
0 & 0 & 0 & 0 & \;0
\end{pmatrix}
\quad \textrm{on } \Sigma_T.
\]
Up to the additional bottom zero row and right zero column this matrix coincides with the corresponding boundary matrix on $\Sigma_T$ in isotropic MHD, cf. \cite{T05MR2187618}.

Then, we easily calculate
\[
\big((\widetilde{\mathcal{A}}_1(\mathring{U},\mathring{\Phi})\dot{V}\cdot\dot{V}\big)\big|_{\Sigma_T}=2\dot{q}(\dot{v}\cdot\mathring{N})|_{\Sigma_T},
\]
where $\dot{q}=\dot{p}_{\bot}+\mathring{H}\cdot \dot{H}$. Similarly to \cite{ST13MR3148595,ST14MR3151094,T10MR2718711}, from the unknowns $\dot{p}_{\bot}$, $\dot{v}$  and $\dot{H}$ we could pass in \eqref{linsym'} to the unknowns $\dot{q}$, $(\dot{v}\cdot\mathring{N},\dot{v}_2,\dot{v}_3)$ and $(\dot{H}\cdot\mathring{N},\dot{H}_2,\dot{H}_3)$ keeping the symmetry of the matrices. But, this is not really necessary. It is only important that the structure of our boundary matrix on the boundary is the same as in isotropic MHD. Since up to the notation $\dot{p}:=\dot{p}_{\bot}$ the boundary conditions \eqref{ELP1c}, \eqref{B'e:def} coincide with those in isotropic MHD, the rest arguments are entirely the same as in \cite{ST13MR3148595,ST14MR3151094}. This completes the proof of Theorem \ref{thm1}.

{\footnotesize 
  }

\end{document}